\documentclass[11pt]{article}
\usepackage{amsmath,amsthm}
\usepackage{pstricks,pst-node,pst-tree}
\usepackage{latexsym}
\usepackage{amssymb,amscd}
\usepackage{geometry, graphicx, subfigure}
\usepackage{url}
\numberwithin{equation}{section}
\let \:=\colon
\let \beg=\begin

\let \mb=\mathbb
\let \mc= \mathcal

\let \ra=\rightarrow

\let \Ga=\Gamma
\let \al=\alpha

\let \fl=\flushleft
\let \fr=\frac

\let \ov=\overline
\let \part=\partial

\let \sub=\subset

\beg{document}
\thispagestyle{empty}
\title{Counting maps from curves to projective space via graph theory}
\vspace{-10pt}
\author{Ethan Cotterill, CMUC, FCT grant PTDC/MAT/111332/2009}
\date{\empty}
\maketitle

\vspace{-20pt}
\section{Brill--Noether theory on reducible curves}
In Brill--Noether theory, one studies linear series on curves, in order to understand when a curve $C$ of genus $g$ comes equipped with a nondegenerate morphism of degree $d$ to $\mb{P}^r$. For a {\it general} curve $[C] \in \mc{M}_g$, a basic answer is provided by the {\it Brill--Noether theorem} of Griffiths and Harris, which establishes that $C$ admits such a morphism if and only if the invariant
\[
\rho(d,g,r)= g-(r+1)(g-d+r)
\]
is nonnegative, in which case $\rho$ also computes the dimension of the space of linear series $g^r_d$ of degree $d$ and rank $r$ on $C$. 

The Brill--Noether question also admits natural extensions, obtained by imposing incidence conditions on the images of the linear series in question. Namely, given integers $m \geq d$ and $s \geq d-r$, let $\mu(d,r,s):= d-r(s+1-d+r)$ denote the virtual dimension of space of inclusions
\beg{equation}\label{inclusion_of_series}
g^{s-d+r}_{m-d}+p_1 + \dots+ p_d \hookrightarrow g^s_m
\end{equation}
on a fixed curve. When the curve $C$ in question is smooth, and the $g^s_m$ is a subspace $V \sub H^0(C,L)$ of global sections of a line bundle $L$, such inclusions correspond to $d$-tuples of points $p_1, \dots, p_d \in C$ for which the natural evaluation map
\beg{equation}\label{evaluation_map}
\mbox{ev}: V \ra H^0(C,L/L(-p_1-\dots-p_d))
\end{equation}
satisfies $\mbox{rank}(\mbox{ev})= d-r$. Geometrically, such $d$-tuples determine $d$-secant $(d-r-1)$-planes to the image of the $g^s_m$. In \cite{Co1}, we showed that when $\rho=0$ and $\mu<0$, there are {\it no} inclusions \eqref{inclusion_of_series} on a general curve:
\beg{thm}\label{Theorem 1}
If $\rho=0$ and $\mu<0$, then a general curve $C$ admits no linear series $g^s_m$ with $d$-secant $(d-r-1)$-planes. 
\end{thm}

Our proof of Theorem~\ref{Theorem 1} is a natural generalization of the Brill--Noether proof given in \cite[Ch. 5]{HM} and is based on an analysis of (limit) linear series on certain reducible curves of compact type.

\section{Counting secant planes via graph theory}
An immediate corollary of Theorem~\ref{Theorem 1} is that {\it when $\rho=0$ and $\mu=-1$, curves with linear series $g^s_m$ with $d$-secant $(d-r-1)$-planes determine a divisor in $\mc{M}_g$}. The case $r=1$ is particularly natural: in that case, exceptional secant planes correspond to $d$-tuples of points for which the evaluation maps \eqref{evaluation_map} fail to be {\it surjective}. We show \cite[Thm 2]{Co2}:
\beg{thm}\label{Theorem 2}
The coefficients of the homology classes of secant-plane divisors in $\ov{\mc{M}}_g$, realized as linear combinations of standard generators over $\mb{Q}$, are explicit linear combinations of hypergeometric series of type $_3F_2$. 
\end{thm}

The key ingredient for proving Theorem 2, which is of interest in its own right, is the following auxiliary result \cite[Thm 4]{Co1}:
\beg{thm}\label{Theorem 3}
The generating series for the virtual number $N_d$ of $d$-secant $(d-2)$-planes to a degree-$m$ curve $C$ of genus $g$ in $\mb{P}^{2d-2}$ is
\beg{equation}\label{N_d}
\sum_{d \geq 0} N_d(g,m)= \bigg(\fr{2}{(1+4z)^{1/2}+1} \bigg)^{2g-2-m} \cdot (1+4z)^{\fr{g-1}{2}}. 
\end{equation}
\end{thm}
Two ingredients enter into our proof of Theorem~\ref{Theorem 3}. The first is {\it Porteous' formula}, which computes the homology class of the locus of $d$-secant $(d-2)$-planes as a determinant in the Chern classes of the so-called $d$th {\it tautological} bundle $L^{[d]}$ over the $d$th Cartesian product $C^d$, whose fiber over $(p_1,\dots,p_d) \in C^d$ is $H^0(L/L(-p_1-\dots-p_d))$. The second is a combinatorial analysis of the resulting intersection-theoretic formula, which amounts to a weighted count of subgraphs of the complete graph $K_d$ on $d$ vertices.

\section{Linear series on metric graphs}
In the preceding section, graphs naturally arose in connection with counting (secant planes to) morphisms via the formalism of intersection theory. But graph theory also intervenes in a natural way as a result of degeneration, via the passage from a nodal curve to the {\it dual graph} recording the incidences of its components.

There is a theory of complete linear series on metric graphs with $\mb{R}$-valued edge lengths due to Baker--Norine \cite{BN} and Mikhalkin--Zharkov \cite{MZ}. Concretely, a (complete) linear series $|D|$ on a metric graph $\Ga$ is a configuration $D$ of points in $\Ga$, modulo an equivalence relation defined by piecewise-linear functions. Moreover, there is an explicit combinatorial {\it burning algorithm} due to Dhar for computing the rank of a configuration of points $D \in \mbox{Div}(\Ga)$; see \cite{CDPR}.

{\bf \fl Contrasting examples in genus four.} Figure (a) shows two metric graphs of genus 4 (here, as in the remainder of the article, we assume that all weights on vertices are 0). The top graph $\Ga_1$ pictured is planar, and the 3 circles determine a degree-3 configuration $D_1$ of trivial rank. Indeed, a fire that burns from $p$ will be repelled by the 3 points in the support of $D_1$, which then evolve at equal velocity against the incoming fire. Assuming the planar graph has generic edge lengths, a single point $p_1$ of $D_1$ will arrive at a vertex $v_1$ of $\Ga_1$, at which point a fire burning from $p$ will approach $v_1$ (and $p_1$) from 2 distinct directions and all of $\Ga_1$ will burn. By contrast, the configuration $D_2$ of 3 points on the complete bipartite graph $\Ga_2=K_{3,3}$ evolves in such a way that at at no time will any fire based at any point $p$ approach any point in the support of $D_2$ along two distinct directions. It follows that $r(D_1) \geq 1$, and in fact the rank of $D_1$ is precisely 1.

\beg{figure}\label{Fig 1}
\hspace{-60pt}\subfigure[Evolutions of degree-3 configurations along the planar ``wheel" graph and $K_{3,3}$]{\includegraphics[scale=.4]{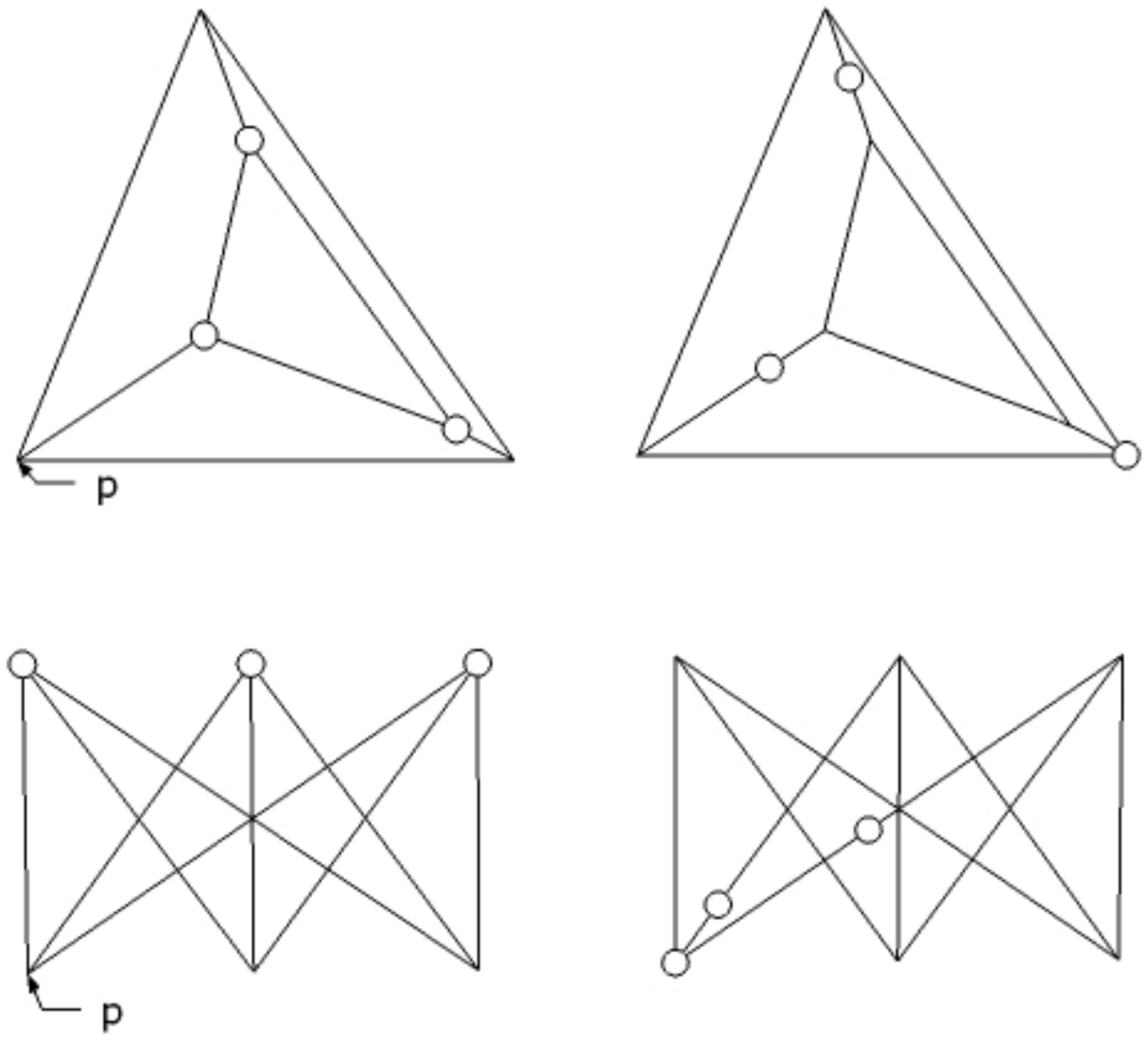}}
\end{figure}

\section{The gonality of tree-decomposed graphs}
The contrast between the behavior of degree-3 configurations on the planar genus-4 graph $\Ga_1$ and on $\Ga_2=K_{3,3}$ is instructive. In fact, it is not hard to check that $\Ga_1$ and $\Ga_2$ each admit two degree-3 configurations of rank 1, as predicted by Brill--Noether theory for curves of genus 4. However, on $\Ga_1$, these configurations depend strongly on the metric structure: each is obtained by placing 2 points on 2 out of 3 inner (resp., outer) ``rim" vertices, and a third point along a ``spoke" at distance from an outer (resp., inner) vertex at distance equal to the length of the shortest spoke. On $\Ga_2$, on the other hand, each rank-1 configuration is associated to 
a choice of one of the two sets of 3 vertices along which $\Ga_2$ decomposes as a union of three 4-edged trees.

{\bf \fl Definition/construction.} Let $V=\{v_1,\dots,v_n\}$, $n \geq 3$ denote a fixed set of vertices, and let $T_1$, $T_2$, and $T_3$ denote three trees each containing $V$ as vertices but which are otherwise pairwise disjoint. The three trees $T_i, 1 \leq i \leq 3$ glue naturally to a graph $\Ga$; we say that $\Ga$ admits a {\it tree decomposition $(T_1,T_2,T_3)$ rooted along $V$}.

Some of the most famous graphs of genus at most 10 admit such tree decompositions: besides $K_{3,3}$, the examples of the so-called Petersen, Heawood, and Pappus graphs in genera 6, 8, and 10 (respectively) are tree-decomposable.

\beg{thm}[Existence of rank-one series on tree-decomposed graphs]\label{existence}
Suppose that the metric graph $\Ga$ admits a tree decomposition rooted on $n \geq 3$ vertices $V$. Then $V$ determines a rank-one, degree-$n$ divisor $D$ on $\Ga$.
\end{thm}

\beg{proof}
The result follows from the burning algorithm. Namely, fix any choice of base point $p$ from which to burn, say $p \in T_1$ without loss of generality. Any fire burning from $p$ along $T_1$ is repelled by the points $p_1, \dots, p_n$ of $D$ supported along $V$, which then evolve at equal velocity along $T_1$ away from $V$. The burning process iterates until ultimately the fire is extinguished by at least one of the points $p_i$, which proves that $r(D) \geq 1$. Similarly, to prove that $r(D)<2$, it suffices to allow {\it two} successive  fires to burn from $p_1$: the first fire simply has the effect of canceling out $p_1$, while the second burns through all of $\Ga$.
\end{proof}

{\bf \fl Definition.} A graph (or a curve) $\Ga$ of genus $g$ is {\it $n$-gonal} whenever $n= \min \{j \in \mb{Z}_{>0}: \exists \text{ a }g^1_j \text{ on }\Ga\}$.

\beg{thm}
$K_{3,3}$, Petersen, Heawood, and Pappus are 3-gonal, 4-gonal, 5-gonal, and 6-gonal graphs, respectively.
\end{thm}

\beg{proof}[Proof sketch]
It is easy to exhibit tree decompositions of these graphs rooted on $n=3, 4$, 5, and 6 vertices, respectively. Whence, by Theorem~\ref{existence}, it suffices to prove that each of these $n$-rooted tree-decomposed graphs admits no degree-$(n-1)$ configurations $D$ of positive rank. Replacing $D$ by a linearly equivalent configuration if necessary, we may assume that each point in $\mbox{Supp}(D)$ appears with multiplicity at most 2. It remains to carry out a case-by-case inspection using the burning algorithm.
\end{proof}

It is not hard to produce graphs that decompose as unions of trees rooted on $n \geq 3$ vertices but are $\al$-gonal with $\al<n$. So additional conditions are needed to ensure that $n$-gonality is achieved. Theorem 4.2 and experimentation give some evidence that it suffices to maximize the minimal cycle length, or {\it girth}, of $\Ga$.

\beg{conj}\label{gonality_conjecture}
A metric graph $\Ga$ that admits a tree-decomposition $(T_1,T_2,T_3)$ rooted on $n$ vertices is $n$-gonal provided $\mbox{girth}(\Ga)$ is {\it maximal} for the combinatorial type of $(T_1,T_2,T_3)$.
\end{conj}

{\bf \fl Acknowledgement.} I am grateful for many illuminating conversations with S. Backman, J. Neves, M. Melo, D. Pinto, and F. Viviani related to linear series on metric graphs.

\beg{thebibliography}{25}
\bibitem{BN} M. Baker and S. Norine, {\it Riemann--Roch and Abel--Jacobi theory on a finite graph}, Adv. Math. {\bf 215} (2007), no. 2, 766--788.
\bibitem{CDPR} F. Cools, J. Draisma, S. Payne, and E. Robeva, {\it A tropical proof of the Brill--Noether theorem}, Adv. Math. {\bf 230} (2012), 759--776.
\bibitem{Co1} E. Cotterill, {\it Geometry of curves with exceptional secant planes: linear series along the general curve}, Math. Zeit. {\bf 267} (2011), no. 3-4, 549--582.
\bibitem{Co2} E. Cotterill, {\it Effective divisors on $\ov{\mc{M}}_g$ associated to curves with exceptional secant planes}, Manuscripta Math. {\bf 138} (2012), no. 1-2, 171--202.
\bibitem{HM} J. Harris and I. Morrison, ``Moduli of curves", Springer, 1998.
\bibitem{MZ} G. Mikhalkin and I. Zharkov, {\it Tropical curves, their Jacobians and theta functions}, Contemp. Math. {\bf 465} (2007), 203--231.
\end{thebibliography}
\end{document}